\documentclass[12pt]{amsart}

\usepackage{epsfig}

\newcommand{\bfz}{{\mathbb {Z}}}

\newtheorem{thm}{Theorem}
\newtheorem{prop}[thm]{Proposition}
\newtheorem{lem}[thm]{Lemma}

\theoremstyle{remark}
\newtheorem{rem}{Remark}

\theoremstyle{definition}
\newtheorem{defn}{Definition}

\newcommand{\Q}{\mathbb{ Q}}

\newcommand{\Z}{\mathbb{ Z}}

\begin{document}

\title[Signatures of surface bundles]{Commutators,
Lefschetz fibrations and the signatures of surface bundles}

\author{H.~Endo}
\address{Department of
Mathematics, Tokyo Institute of Technology, Oh-okayama, Meguro
152-8551, Tokyo, Japan}
\email{endo@math.titech.ac.jp}

\author{M.~Korkmaz}
\address{Department of Mathematics, Middle East Technical University,
06531 Ankara, Turkey}
\email{korkmaz@arf.math.metu.edu.tr}

\author{D.~Kotschick}
\address{Mathematisches Institut, Universit\"at M\"unchen,
Theresienstr.~39, 80333 M\"unchen, Germany}
\email{dieter@member.ams.org}

\author{B.~Ozbagci}
\address{Department of Mathematics, Michigan State University, E. Lansing,
MI 48824-1027, USA}
\email{bozbagci@math.msu.edu}

\author{A.~Stipsicz}
\address{Department of Analysis, ELTE TTK,
M\'uzeum krt. 6-8, H-1088 Budapest, Hungary}
\email{stipsicz@cs.elte.hu}

\date{\today; MSC 2000: 57M05, 20F65, 57R20}

\thanks{The first author was supported by the {\sl Deutsche
Forschungsgemeinschaft}; the third author
is a member of the {\sl European Differential Geometry Endeavour}
(EDGE), and the fifth author was partially supported by OTKA
and Sz\'echenyi Professzori \"Oszt\"ond{\'\i}j.}

\begin{abstract}
We construct examples of Lefschetz fibrations
with prescribed singular fibers. By taking differences of pairs
of such
fibrations with the same singular fibers, we obtain new
examples of surface bundles over surfaces with non-zero
signature. From these we derive new upper bounds for the
minimal genus of a surface representing a given element
in the second homology of a mapping class group.
\end{abstract}

\maketitle


\newenvironment{prooff}{\medskip \par \noindent {\it Proof}\ }{\hfill
$\square$ \medskip \par}
    \def\sqr#1#2{{\vcenter{\hrule height.#2pt
        \hbox{\vrule width.#2pt height#1pt \kern#1pt
            \vrule width.#2pt}\hrule height.#2pt}}}
    \def\square{\mathchoice\sqr67\sqr67\sqr{2.1}6\sqr{1.5}6}

\section{Introduction}
\label{first}
It is an elementary fact that the Euler characteristic is multiplicative
in fiber bundles. According to a classical result of Chern, Hirzebruch
and Serre~\cite{CHS} the same holds for the signature, provided that the
fundamental group of the base acts trivially on the cohomology of the
fiber. Atiyah~\cite{A} and, independently, Kodaira~\cite{Kd} showed
that this assumption on the monodromy is necessary, by exhibiting
surface bundles over surfaces with non-zero signature.

In the case of bundles whose fiber is a sphere or torus, it is
easy to see that the signature must vanish. Therefore, only the
signature of surface bundles of higher genus is interesting.
For a closed oriented surface $F$ of genus $h\geq 2$, Teichm\"uller
theory implies that the identity component
of the group of orientation-preserving diffeomorphisms is
contractible. It follows that every oriented bundle with fiber
$F$ over a base $B$ is determined by (the conjugacy class of)
its monodromy representation
$$
\rho\colon\pi_1(B)\longrightarrow\Gamma_h \ ,
$$
where $\Gamma_h$ is the mapping class group of $F$, consisting
of isotopy classes of orientation-preserving diffeomorphisms.
If the base $B$ is also $2$-dimensional, then the signature
of the total space $X$ is four times the first Chern number
of the flat symplectic bundle obtained by composing $\rho$
with the action of $\Gamma_h$ on the homology of $F$,
see~\cite{A,Ho}. In particular, the signature vanishes if
the genus of $B$ is $0$ or $1$.
The signature also vanishes for all bundles with fiber genus
$2$, because of Igusa's theorem $H_2(\Gamma_2,\Q )=0$.
Thus, we may assume that the fiber genus $h$ is $\geq 3$.

Combining the work of Meyer~\cite{M} and of Harer~\cite{H}, one
sees that the signature of the total space $X$ is given by the
homology class of $\rho_*[B]$ in the homology of $\Gamma_h$.
More precisely, the
second integral homology of the mapping class
group is infinite cyclic, generated by the Meyer signature cocycle
corresponding to the signature of the total space.
This means that determining the maximal signature of a surface
bundle with given fiber and base genus is equivalent to calculating
the Gromov-Thurston norm in the second homology of the mapping
class group. This is essentially Problem 2.18 in Kirby's list~\cite{Ki}.
To address this problem, consider the function
$$
g_h (n)=\min \{ g \mid \exists \ \mbox{ a $\Sigma _h$-bundle }
X\to \Sigma _g \mbox{ with } \sigma (X)=4n\} \ .
$$

Using Seiberg-Witten gauge theory, the first nontrivial lower bound
for this function was proved in~\cite{Kt}:
\begin{equation}\label{Kot}
g_h(n)\geq \frac{2\vert n\vert }{h-1} +1 \ .
\end{equation}
The only systematic upper bound for this function was proved in~\cite{E},
where it was shown that for every fiber genus $h\geq 3$ there is
a surface bundle over a surface of genus $111$ with signature $4$.
Pulling back to coverings of the base, one has
\begin{equation}\label{Endo}
g_h(n)\leq 110\vert n\vert +1 \ .
\end{equation}
A non-explicit improvement of~\eqref{Endo} in some cases was
given in~\cite{S}.

In this paper we obtain new upper bounds for the function $g_h(n)$
by constructing examples of surface bundles in which the base genus
is comparatively small. We found these examples by
first constructing Lefschetz fibrations with singular fibers
corresponding to expressions of products of Dehn twists as products
of commutators, and then taking differences of Lefschetz fibrations
with the same singular fibers to obtain smooth surface bundles.
We have chosen to present the examples in the way we originally
found them, although it would have been possible, after the fact, to
eliminate the Lefschetz fibrations from the presentation and write
down the monodromy representations of the surface bundles directly.
We believe that the subtraction of Lefschetz fibrations presented
in Section~\ref{second}, also used in~\cite{S}, is of interest in its
own right, in addition to being a useful stepping stone in the
construction of surface bundles.

Our first main theorem is the following improvement of~\eqref{Endo}:
\begin{thm}\label{main1}
For every $h\geq 3$ there is a surface bundle of genus $h$
over the surface of genus $9$ with signature $4$.
In particular, $g_h (n)\leq 8\vert n\vert +1$.
\end{thm}

Notice that all these examples over $\Sigma _{9}$ have the same signature.
By considering sections of our fibrations we can construct surface
bundles with fiber genus $h$ over $\Sigma _{9}$ for which the signature grows linearly with $h$. More precisely, we have:
\begin{thm}\label{main2}
For every $h\geq 3 $ there are surface bundles of fiber genus $h$
over the surface of genus $9$ with signature at least
$4\frac{h-2}{3}$.
\end{thm}
This result allows us to prove upper bounds for $g_h(n)$ which
have the same shape as the lower bound~\eqref{Kot}, in that
the fiber genus appears in the denominator. We only formulate
these upper bounds in the asymptotic case, when $n$ becomes
large.
It is easy to see that the limit
$$
G_h=\lim _{n\to \infty } \frac{g_h (n)}{n}
$$
exists and is finite for all $h$. The inequality~\eqref{Kot}
implies $G_h \geq \frac {2}{h-1}$.
Using our new examples, we will prove:
\begin{thm}\label{main3}
If $h\geq 3$ is odd, then $G_h\leq \frac{16}{h-1}$.
If $h\geq 4$ is even, then $G_h\leq \frac{16}{h-2}$.
\end{thm}

This paper is organized as follows.
In Section~\ref{second} we review the basic facts about Lefschetz
fibrations and
describe the ``subtraction operation" for them in detail. Section~\ref{third}
is devoted to the proof of various identities in the mapping class group
expressing certain products of Dehn twists as products of commutators.
In Section~\ref{fourth} we calculate the signatures of the
corresponding Lefschetz fibrations using the Meyer signature cocycle~\cite{M}.
In the last Section we give the proofs of the Theorems stated above.

\section{Subtracting Lefschetz fibrations}
\label{second}
We begin by recalling the definition and basic properties of
Lefschetz fibrations. More details can be found in~\cite{GS,La}.
Let $X$ be a compact oriented $4$-manifold, and $B$ a compact
oriented surface.
\begin{defn}
A smooth map $f\colon X\to B$ is called a Lefschetz fibration
if it is surjective and if for each critical point $p \in X$
there are local complex coordinates $(z_1,z_2)$ on $X$ around $p$
and $z$ on $B$ around $f(p)$ compatible with the orientations
and such that $f(z_1,z_2)=z_1^2+z_2^2$.
\end{defn}
It follows that a Lefschetz fibration
has at most finitely many critical points $p_1, \ldots , p_k$.
It is easy to see that by a slight perturbation one can achieve
that $f$ is injective on its critical set
$C = \{ p_1, \ldots, p_k \}$. We will always assume that this
additional property holds.

The {\em genus} of $f$ is defined to be the genus of a regular fiber.
If $B$ is connected, the genus is well-defined. Even when $B$ is
not connected, we will assume that all regular fibers have the
same genus. Fibers of $f$ passing through elements of $C$ are
{\em singular} fibers. Notice that if $\nu (f(C))$ denotes an open
tubular neighborhood of the set of critical values $f(C)$,
then the restriction of $f$ to $f^{-1}(B\setminus\nu (f(C)))$
is a smooth surface bundle over the surface-with-boundary
$B\setminus \nu (f(C))$.

A singular fiber $f^{-1}(q_i)$, where $q_i=f(p_i)$, can be
described by its {\em monodromy}, which is an element in the
mapping class group $\Gamma_h$. To determine this element,
however, we need to fix a base point $\tau \in B\setminus f(C)$,
an identification  of $f^{-1}(\tau )$ with the closed oriented
surface $F$ of genus $h$, and a loop $c_i$ in $B$  based at
$\tau$ which has linking number $+1$ with $q_i$. The restriction
of $f$ to the preimage of this loop is an $F$-bundle over $S^1$
which can be described by a single
element $t_i\in \Gamma_h$. In fact, by performing this
procedure for all loops in $B\setminus f(C)$ we get a map
$\varphi \colon \pi _1 (B\setminus f(C))\to \Gamma_h$.
It can be shown that $t_i$ is a right-handed Dehn twist
along a simple closed curve $v_i \subset f^{-1}(\tau )$
called the {\em vanishing cycle} corresponding to the singular
fiber $f^{-1}(q_i)$. Notice that, even after fixing
$\tau \in B$ and the identification $F \approx f^{-1}(\tau )$,
both $t_i$ and $v_i$ depend on the chosen loop $c_i$.

\begin{figure}[hpbt]
\centerline{\epsfbox{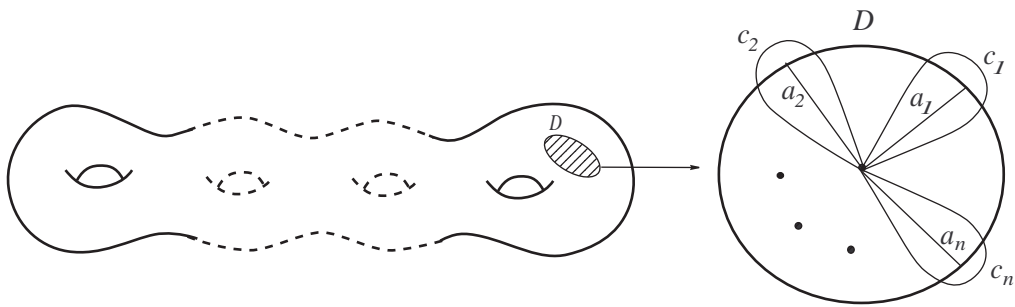}}
\caption{Choice for loops defining vanishing
cycles.} \label{loops}
\end{figure}

It is convenient to fix the following conventions. Suppose that all
$q_i$ lie on the boundary of a  disk $D\subset B$ centered at
$\tau \in B$. Let $a_i$ denote the radial curve in $D$ connecting $\tau$
with $q_i$ and form $c_i$ as the boundary of an appropriate neighborhood
of $a_i$, cf.~Figure~\ref{loops}.
By fixing a generating system $\{ a_1, b_1,\ldots , a_g , b_g\}$ of
$\pi _1 (B\setminus D)$, the map  $\varphi$ can be encoded
by a sequence $(t_1, \ldots , t_s, \alpha _1, \beta _1, \ldots ,
\alpha _g , \beta _g )$, where $\alpha _i $ and $\beta _i\in \Gamma_h$
tell us the monodromy of the fibration along $a_i $ and $b_i$.
It is easy to see that these elements satisfy the relation
$\Pi _{j=1}^g [\alpha _i , \beta _i]\cdot\Pi _{i=1}^s  t_i =1$ in
the mapping class group. Conversely, for $h\geq 2$ a word of the form
$\Pi _{j=1}^g [\alpha _i , \beta _i]\cdot\Pi _{i=1}^s t_{i}$ representing
1 in $\Gamma_h$ (with  $t_i$ being right-handed Dehn twists)
gives rise to a Lefschetz fibration of genus $h$ over a
surface $B$ of genus $g$.

\begin{figure}[hpbt]
\centerline{\epsfbox{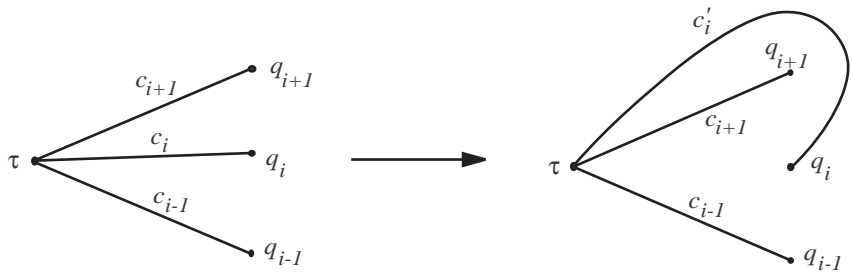}}
\caption{Elementary transformation.}
\label{trafo}
\end{figure}

As we noted already, the vanishing cycles and the corresponding Dehn
twists depend on the chosen loops $c_i$. It is easy to see that
a cyclic permutation of the indices can be compensated by changing the
identification $F\approx f^{-1}(\tau )$, so the resulting Lefschetz
fibration remains the same. One can also change the word by
{\em elementary transformations} without changing the Lefschetz
fibration, i.~e.~the path $c_i$ can be changed as indicated
by Figure~\ref{trafo}.
By applying an elementary transformation as shown by the figure, we replace
$t_i$ and $t_{i+1}$ by $t_{i+1}$ and $t_{i+1}^{-1}t_i t_{i+1}$.
(Notice that this change has no effect on the product of these elements.)
The new vanishing cycles are easy to determine since for any mapping class
$g$ the conjugate $g^{-1}t_vg$ of the Dehn twist $t_v$ is simply the
Dehn twist $t_{g(v)}$. It is not hard to prove that if two words
give rise to equivalent fibrations then the words can be transformed
into each other by applying combinations of the two operations just
described.

A singular fiber $f^{-1}(q_i)$ is {\em nonseparating} if the corresponding
vanishing cycle $v_i$ is nonseparating, equivalently its homology class is
nonzero in $H_1(f^{-1}(\tau ); \bfz )$. If $v_i$ is a separating curve,
equivalently its homology class is zero, then $f^{-1}(q_i)$ is called
{\em separating}. A vanishing cycle $v_i$ and the corresponding singular
fiber are of {\em type} $0$ if $v_i$ is nonseparating; they are of {\em type}
$j\in \{1,\ldots , [\frac{h}{2}]\} $ if the vanishing cycle separates the
surface of genus $h$ into two components with genera $j$ and $h-j$. Although
the vanishing cycle depends on the chosen path $c_i$, its type is independent
of this choice. From the classification of surfaces, one can prove that for
two simple closed curves of the same type there exists a diffeomorphism
of the ambient surface mapping one into the other. This implies that two
singular fibers of the same type have  fiber- and orientation-preservingly
diffeomorphic tubular neighborhoods.

The combinatorial data of a Lefschetz fibration can be encoded as follows:
\begin{defn}\label{d:comb}
The vector $\mu _{comb}(X)= (\mu _0, \ldots , \mu _{[\frac{h}{2}]})\in
\bfz ^{[\frac{h}{2}]+1}$ associated to the Lefschetz fibration
$f\colon X\to B$ is constructed by taking $\mu _j$ to be the number
of singular fibers of type $j$ ($j=0, 1, \ldots , [\frac{h}{2}]$).
Following~\cite{Sm} we say that two fibrations
$f_i \colon X_i\to B_i$ ($i=1,2$) are combinatorially equivalent
if $\mu _{comb}(X_1)=\mu _{comb}(X_2)$.
\end{defn}

The construction we use to produce new examples of surface
bundles is a procedure for taking the difference of two combinatorially
equivalent Lefschetz fibrations. If $X_1$ and $X_2$ are
combinatorially equivalent as in Definition~\ref{d:comb},
with critical values $\{ q^1 _i\}_{i=1}^s$ and $\{ q^2 _i\}_{i=1}^s$
respectively, then a surface bundle $X_1-X_2\to B_1-B_2$ can be constructed
in the following way: order the $q^1_i$'s and $q^2_j$'s so that
singular fibers with coinciding lower index have the same type.
Fix an orientation- and fiber-preserving diffeomorphism $\phi _i$ between the
boundaries of tubular neighborhoods of fibers with lower index $i$
($i=1,\ldots , s$). The union of these maps will be denoted by $\phi$.
Now glue $X_1\setminus (\cup_{i=1}^s \nu (f_1(q^1 _i)))$ to
${\overline{X_2\setminus (\cup_{i=1}^s \nu (f_2(q^2 _i)))}}$
using $\phi$. Notice that by reversing the orientation on $X_2$,
the map $\phi$ becomes orientation-reversing, hence the resulting
manifold $Y=X_1-X_2$ inherits a natural orientation. Since $\phi $
is fiber-preserving, $Y$ admits a smooth fibration with fibers of
genus $h$ over a compact surface $B$ which we will denote by $B_1-B_2$.

\begin{lem}\label{l:sub}
If $X_1\to B_1$ and $X_2\to B_2$ are combinatorially equivalent
Lefschetz fibrations with $s$ singular fibers,
then $Y=X_1-X_2$ is a smooth surface bundle with signature
$\sigma (Y)=\sigma (X_1)-\sigma (X_2)$ over the surface $B = B_1-B_2$
with Euler characteristic $\chi (B) =\chi (B_1)+\chi (B_2)-2s$.
\end{lem}
\begin{proof}
By construction, $Y$ is an oriented smooth surface bundle over a
surface $B$. The claim about the Euler characteristic of the base
is obvious. The claim about
the signature is an instance of Novikov additivity.
\end{proof}

Note that we did not assume $X_1$ and $X_2$ to be connected. This
means that basepoints have to be chosen in each component of $B_i$,
and the vectors of combinatorial data have to be summed over all
components to determine combinatorial equivalence.
If $X_1$ or $X_2$ happens to be connected, then so is $X_1-X_2$.

The main property we used in the above construction is that
the manifolds $X_1\setminus (\bigcup _{i=1}^s \nu (f_1(q^1 _i)))$
and $X_2\setminus (\bigcup _{i=1}^s \nu (f_2(q^2 _i)))$ have
diffeomorphic boundaries and, after reversing the orientation of
one of them, this diffeomorphism can be chosen to be fiber-preserving
and orientation-reversing.
A variation of this construction goes as follows:
Suppose that  partitions of the critical values $\{ q^1_i\}_{i=1}^s$ and
$\{ q^2_i\}_{i=1}^s$ are given together with a system of disjoint disks $D^k _j
\subset B_k$ ($k=1,2$ and $j=1, \ldots ,m$)
such that each disk contains exactly one equivalence
class of the partitions.
Suppose furthermore that we can pair up these disks
in a way that the surface bundles
$X_1\vert _{ D^1_j}$ are isomorphic to $X_2 \vert _{ D^2_j}$
for all $j=1,\ldots , m$. Then $X_2$ can be subtracted
from $X_1$ along the disks $D^k_j$, i.~e.~the manifold
$$
Y=(X_1\setminus (\bigcup_{j=1}^m f_1 ^{-1} ({\rm {int }}D^1_j)))\bigcup
{\overline{(X_2\setminus (\bigcup_{j=1}^m f_2 ^{-1} ({\rm {int }}D^2_j)))}}
$$
admits the structure of a  surface bundle. The signature $\sigma (Y)$
is again given by $\sigma (X_1)-\sigma (X_2)$, while the Euler
characteristic of the base is equal to $\chi (B_1)+\chi (B_2)-2m$.

\begin{rem}
The definition of $X_1-X_2$ is a special case of this latter construction,
corresponding to the situation when each equivalence class of the
partition consists of a unique critical value. By considering a partition
with larger equivalence classes we get  smaller $m$ which results in a smaller
genus for the base. Notice that in the special case of $X_1-X_2$ the
assumption $X_1\vert _{ D^1_j} \approx X_2\vert _{ D^2_j}$
can be  easily checked by determining the type of the singular fibers
over the disks. In general, however, the types of the singular
fibers over the disks do not specify the diffeomorphism
type of the above fibration, since
fibers of the same type can be glued together in many different ways
resulting various fibrations over $ D_j ^k$.
\end{rem}

\begin{rem}
There is a generalisation of Lefschetz fibrations, called achiral
Lefschetz fibrations, where one allows singular fibers whose
monodromies are left-handed Dehn twists, cf.~\cite{GS}. Keeping
track of the chirality of the singular fibers, it is clear that
the subtraction operation described above generalises to the
category of achiral Lefschetz fibrations.
\end{rem}

We conclude this section by discussing the relation between the word
specifying a Lefschetz fibration and sections of the fibration.
Suppose that $f\colon X\to B $ is a given  Lefschetz
fibration. A map $\sigma \colon B \to X$ is called a {\em section}
if $f\circ \sigma ={\rm {id}}_{B }$. The {\em self-intersection}
 (or {\em square}) of the section $\sigma$ is simply the self-intersection
number of the homology class $[\sigma (B )]\in H_2(X; \bfz )$.
In the following $\Gamma_{h,1}$ denotes the mapping class group of the
closed oriented surface of genus $h$ with one marked point
and $\Gamma^1_h$ denotes the mapping class group with respect to one
boundary component (fixed pointwise). Notice that by collapsing the
boundary circle to a point we get a natural surjection
$\varphi \colon \Gamma^1_h\to \Gamma_{h,1}$ with kernel the subgroup
generated by the Dehn twist $\Delta _{\partial}$
along a curve isotopic to the boundary circle (cf.~\cite{Waj}, for
example). Moreover, by forgetting the marked point we have an obvious map
$\Gamma_{h,1}\to \Gamma_h$.

The following two well-known facts show how the existence of a section
(and its square) is reflected in the monodromy representation of a
Lefschetz fibration.
Suppose that the monodromy representation of $f\colon X\to B$
is given by the relator $\Pi _{j=1}^g [a_i, b_i]\cdot\Pi _{i=1}^st_i$
representing 1 in $\Gamma_h$.

\begin{prop}
The fibration admits a section if and only if $t_i$ and $a_j, b_j\in
\Gamma_h$ admit lifts
${\tilde {t}}_i, {\tilde {a}}_j, {\tilde {b}}_j \in \Gamma_{h,1}$
such that
$\Pi _{j=1}^g [{\tilde {a}}_j, {\tilde {b}}_j] \cdot \Pi _{i=1}^s {\tilde {t}}_i$ represents $1$ in $\Gamma_{h,1}$. A section of $f\colon X\to B $
is given once such a lift is fixed. \qed
\end{prop}

Suppose now that a fibration $f\colon X\to B$ with a section is given,
so a  lift $\Pi _{j=1}^g [{\tilde {a}}_j, {\tilde {b}}_j]\cdot
\Pi _{i=1}^s {\tilde {t}}_i$ of $\Pi _{j=1}^g [a_j, b_j]\cdot
\Pi _{i=1}^s t_i $ is fixed. Take a lift $t'_i$ of ${\tilde {t}}_i$
(and $a_j', b_j'$ of ${\tilde {a}}_j, {\tilde {b}}_j$) in $\Gamma^1_h$
and consider $\Pi _{j=1}^g [a_j' , b_j']\cdot\Pi _{i=1}^s t_i '
\in \Gamma_h ^1$. By the discussion
above, this product is in $\ker \varphi$, hence it is equal to
$\Delta ^n _{\partial}$ for some $n \in \bfz $.

\begin{prop}{\rm  {(cf.~\cite{Sm})}}
The self-intersection number of the section given by the
above lift is equal to $-n$. \qed
\label{secc}
\end{prop}
Next we would like to show that after subtracting Lefschetz
fibrations with sections, under favourable circumstances
the resulting fibration admits a section whose self-intersection
number is equal to the difference of the self-intersection
numbers of the sections of the individual fibrations.
For this, suppose that two fibrations $f_i \colon X_i \to B _i $ ($i=1,2$)
are given by their monodromy representations
$\Pi [a_j, b_j]\cdot\Pi t_i $ and $\Pi [c_j , d_j]\cdot\Pi s_j $
respectively. Suppose
furthermore that the disks $D_i\subset B _i$ along which the subtraction
operation  will be performed contain the singular fibers corresponding to
the Dehn twists $t_{i_1}\ldots t_{i_k}$ (and   $s_{i_1}\ldots s_{i_k}$ resp.).

\begin{prop}\label{osszerag}
If the lifts ${\tilde {t}}_{i_n}$ giving rise to the sections  coincide
with ${\tilde {s}}_{i_n}$ ($n=1,\ldots , k$) in $\Gamma_{h,1}$,
then the difference of the two fibrations admits a section.
The self-intersection of this section is given
by the difference of the self-intersection of the individual pieces.
\end{prop}
\begin{proof}
The assumption shows that there is a diffeomorphism $f^{-1}_1(D_1)\to
f^{-1}_2(D_2)$ mapping  the sections  into each other. Now the statement
is obvious --- notice that the self-intersection is the difference of
the two self-intersections since in the subtracting operation we
change the orientation of $X_2$.
\end{proof}
\begin{rem}
The assumption on coinciding lifts cannot be relaxed, as the following
example shows. Take two copies of the trivial bundle
$\Sigma \times S^2 \to S^2$, fix two sections in each and blow up one section
in each copy. In this way we get two Lefschetz fibrations (each with a
single singular fiber) for which the subtraction operation (along the
singular fibers) applies and gives $\Sigma _h \times S^2 \to S^2$ back.
The section blown up, however, can be glued only to the section in the
other copy also blown up, because otherwise we would find a homology
class in
$\Sigma _h \times S^2$ with odd square, which is clearly impossible.
\end{rem}
Surface bundles with sections of self-intersection zero can be summed
along their sections by performing a fiberwise connected sum. This
is an instance of Gompf's symplectic sum operation, but for our
purposes the symplectic aspect is irrelevant.
\begin{lem}\label{secrag}
If $X_i\to B$ with $i=1,2$ are two surface bundles with fiber genera
$h_i$ over the same base surface and both fibrations admit sections
with self-intersection zero, then there is a
surface bundle over $B $ with fiber genus $h_1+h_2$ and
signature $\sigma (X_1)+\sigma (X_2)$.
\end{lem}
\begin{proof}
The signature is additive when summing along embedded surfaces of
self-intersection zero.
\end{proof}

\section{Commutators in mapping class groups}
\label{third}
Let $F_{h,s}^r$ be an oriented surface of genus $h$ with $s$
marked points and $r$ boundary components. The {\em mapping class
group} $\Gamma_{h,s}^r$ of $F$ consists of the isotopy classes
of orientation-preserving diffeomorphisms of $F$ which are the
identity on each boundary component and preserve the set of
marked points. The isotopies are not allowed to permute marked
points or to rotate boundary components. The groups $\Gamma_{h,s}^0$,
$\Gamma_{h,0}^r$ and $\Gamma_{h,0}^0$ will be denoted by
$\Gamma_{h,s}$, $\Gamma_{h}^r$ and $\Gamma_h$, respectively.

We say that two simple closed curves $a$ and $b$ on $F$ are
{\em topologically equivalent} if there exists a diffeomorphism of $F$
mapping $a$ to $b$. For a group $G$ and $x,y\in G$, the commutator
$[x,y]$ denotes the element $xyx^{-1}y^{-1}$ and $x^y$ denotes the
conjugate $yxy^{-1}$.

It follows easily from the definition of a Dehn twist that if $a$ is a
simple closed curve on $F$ and $f$ is an orientation-preserving
diffeomorphism of $F$, then $ft_af^{-1}=t_{f(a)}$ in $\Gamma_{h,s}^r$.

If $a$ and $b$ are two topologically equivalent simple closed curves
on $F$, then $t_at_b^{-1}$ is a commutator. More precisely, if $f(a)=b$
then $t_at_b^{-1}=[t_a,f]$.

Let $a$ and $b$ be two simple closed curves on $F$. If $a$ is disjoint
from $b$, then the supports of the Dehn twists $t_a$ and $t_b$ can be
chosen to be disjoint. Hence, $t_a$ commutes with $t_b$. If $a$ intersects
$b$ transversely at one point, then it is easy to see that $t_at_b(a)=b$.
It follows that $t_a$ and $t_b$ satisfy the {\em braid relation}
$t_at_bt_a=t_bt_at_b$.

\begin{figure}[hbt]
    \begin{center}
          \epsfig{file=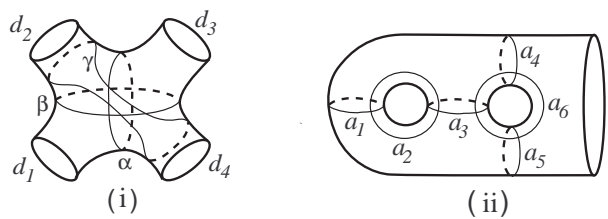,width=8.5cm}
          \caption{Curves of the lantern and two-holed torus relation.}
          \label{lantern}
    \end{center}
\end{figure}

The following two relations in the mapping class group are also well-known.
The first one is the {\em lantern relation} (cf.~\cite{J}). Let $S$ be a
sphere with four boundary components $d_1$, $d_2$, $d_3$ and $d_4$. Suppose
that $S$ is embedded in $F$. Then there are three simple closed curves
$\alpha$, $\beta$, $\gamma$ on $S$, as illustrated in Figure~\ref{lantern}~(i),
which satisfy the lantern relation
$$
t_{d_1}t_{d_2}t_{d_3}t_{d_4}=t_{\alpha}t_{\beta}t_{\gamma} \ .
$$

The second relation is the {\em two-holed torus relation} or {\em chain
relation}. Let $a_1$, $a_2$, $a_3$ be three nonseparating simple closed curves
on $F$ such that $a_2$ intersects $a_1$ and $a_3$ transversely only once,
$a_1$ is disjoint from $a_3$ and $a_1\cup a_3$ does not disconnect $F$. A
regular neighbourhood of $a_1\cup a_2\cup a_3$ is a torus with two nonseparating
boundary components, say $a_4$ and $a_5$ (cf.~Figure~\ref{lantern}~(ii)).
Clearly, $a_4$ and $a_5$ are disjoint from $a_1$, $a_2$, $a_3$, and from each
other. By using the braid relation and the fact that $t_{a_1}$ commutes with
$t_{a_3}$, the relation given by Proposition~3 in~\cite{Li} is easily shown
to be equivalent to the two-holed torus relation
$$
t_{a_4}t_{a_5}=(t_{a_1}t_{a_2}t_{a_3})^4 \ .
$$

Now we describe various commutator relations in mapping class groups.
These relations will be used in the next Section to construct the Lefschetz
fibrations used in the course of the proofs of the Theorems stated
in Section~\ref{first}.


\begin{lem}\label{Proposition}
Let $a$, $b$, $c$ and $d$ be four simple closed curves on $F$ such that
$a$ is disjoint from $b$, $c$ is disjoint from $d$, and the complements
of $a\cup b$ and $c\cup d$ in $F$ are connected.
Then $t_at_b^{-1} t_ct_d^{-1}$ is a commutator.
\end{lem}
\begin{proof}
By the classification of surfaces, there exists a diffeomorphism $g$
of $F$ such that $g(a)=d$ and $g(b)=c$. Then
$$
t_at_b^{-1} t_ct_d^{-1}=t_at_b^{-1} t_{g(b)}t_{g(a)}^{-1}
=t_at_b^{-1} gt_bt_a^{-1} g^{-1}
=[t_at_b^{-1}, g] \ .
$$
\end{proof}


\begin{prop}\label{thm1}
Let $h\geq 3$ and let $a$
be a simple closed curve on $F$. In
 the mapping class group $\Gamma_{h,s}^r$ of $F$
 \begin{itemize}
  \item[(a)] $t_a^2$ can be written as a product of two commutators,
  \item[(b)] if $a$ is nonseparating, then $t_a^4$ can be written as
  a product of three commutators.
 \end{itemize}
\end{prop}
\begin{proof}
Suppose that the surface of genus $3$ with two holes in Figure~\ref{feners}
is embedded in $F$.
\begin{figure}[hbt]
    \begin{center}
          \epsfig{file=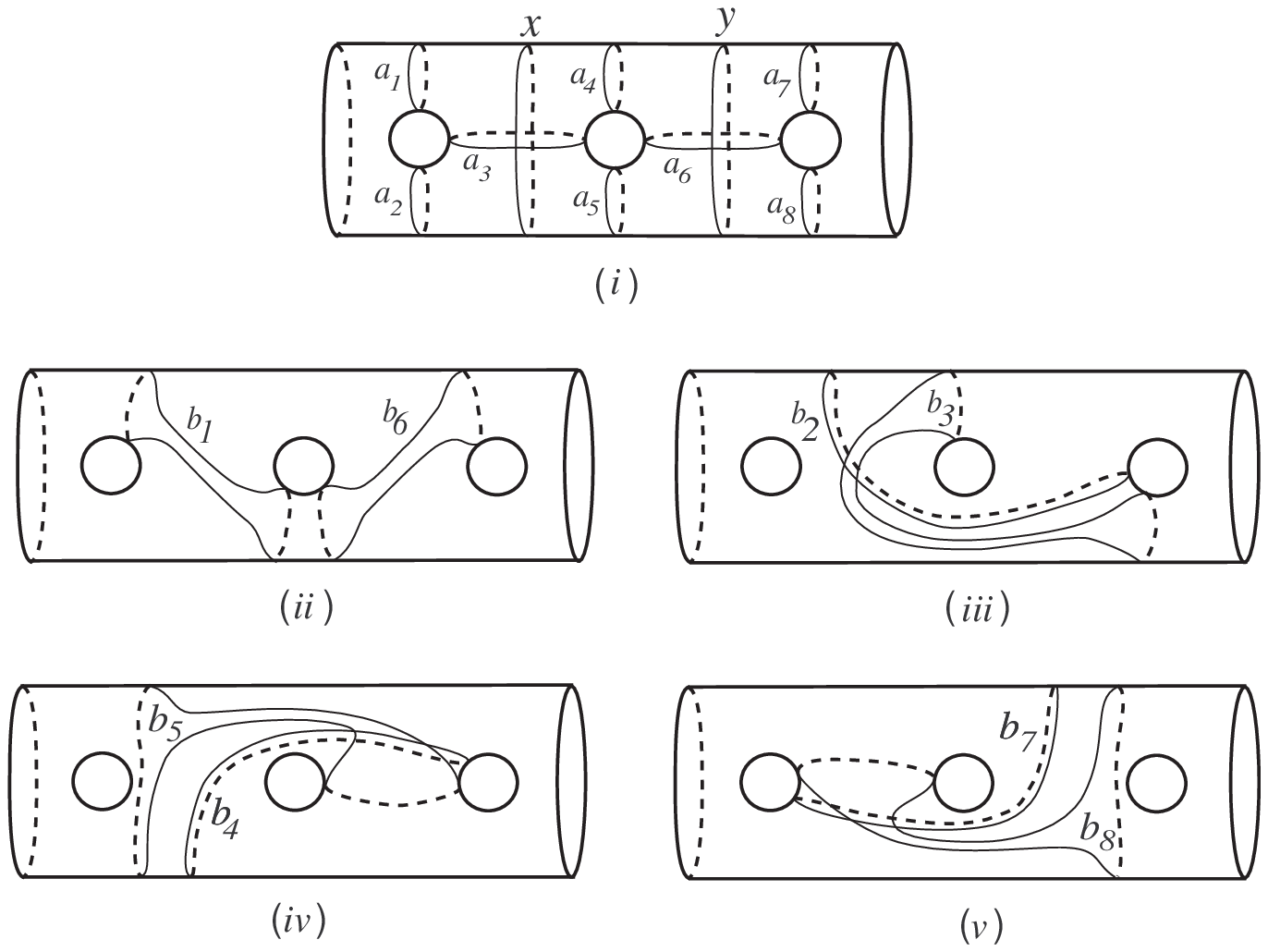,width=12cm}
          \caption{}
          \label{feners}
    \end{center}
 \end{figure}
Consider the curves on $F$ given in the figure. The sphere $S$ with four
holes of the lantern relation, see Figure~\ref{lantern}, can be embedded
in $F$ so that the curves
$d_1$, $d_2$, $d_3$, $d_4$, $\alpha$, $\beta$, $\gamma$ become
respectively $a_2$, $a_1$, $a_4$, $a_5$, $x$, $a_3$, $b_1$. This gives
us the relation
\begin{equation}\label{eqn3}
   t_{a_1}t_{a_2}t_{a_4}t_{a_5}=t_{x}t_{a_3}t_{b_1} \ .
\end{equation}
Similarly, two other embeddings of $S$ give the relations
\begin{equation}\label{eqn4}
   t_{x}t_{a_4}t_{a_6}t_{a_8}=t_{a_5}t_{b_2}t_{b_3}
\end{equation}
 and
\begin{equation}\label{eqn5}
   t_{x}t_{a_5}t_{a_6}t_{a_7}=t_{a_4}t_{b_4}t_{b_5} \ .
\end{equation}
If we multiply both sides of~\eqref{eqn3} by $t_{b_2}t_{b_3}$,
use~\eqref{eqn4} and cancel $t_x$, we obtain
\begin{equation}\label{*}
t_{a_4}^2t_{a_1}t_{a_2}t_{a_6}t_{a_8}=t_{a_3}t_{b_1}t_{b_2}t_{b_3} \ ,
\end{equation}
or, equivalently,
\begin{equation}\label{(a1)}
t_{a_4}^2=t_{a_3}t_{a_6}^{-1}\,  t_{b_1}t_{a_8}^{-1}\,
    t_{b_2}t_{a_1}^{-1}\,  t_{b_3}t_{a_2}^{-1}.
\end{equation}
Similarly, the equalities~\eqref{eqn4} and~\eqref{eqn5} yield the
equality
\begin{equation}\label{(a2)}
t_{x}^2=t_{b_2}t_{a_6}^{-1}\,  t_{b_3}t_{a_7}^{-1}\,
    t_{b_4}t_{a_6}^{-1}\,  t_{b_5}t_{a_8}^{-1} \ .
\end{equation}
Applying Lemma~\ref{Proposition} to~\eqref{(a1)} and~\eqref{(a2)}
proves that $t_{a_4}^2$ and $t_{x}^2$ are products of two commutators.
Any nonseparating simple closed curve is topologically equivalent to
$a_4$. If $a$ is a separating simple closed curve on $F$, then the
surface of genus $2$ on the right hand side of $x$ can be embedded in
$F$ so that $x$ is topologically equivalent to $a$. Now, the proof of
(a) follows from the fact that a conjugate of a commutator is again a
commutator.

Similarly, two more embeddings of the lantern give the relations
\begin{eqnarray}
    t_{a_4}t_{a_5}t_{a_7}t_{a_8}=t_{b_6}t_{a_6}t_{y}
      \label{eqn6}\\
    t_{y}t_{a_2}t_{a_3}t_{a_4}=t_{a_5}t_{b_7}t_{b_8}
      \label{eqn7}.
\end{eqnarray}
Multiplying~\eqref{eqn6} by $t_{b_7}t_{b_8}$ from the left and
using~\eqref{eqn7} gives
\begin{eqnarray}
    t_{a_4}^2 t_{a_2}t_{a_3}t_{a_7}t_{a_8}
    =t_{b_7}t_{b_8}t_{b_6}t_{a_6}.  \label{eqn8}
\end{eqnarray}
By combining~\eqref{*} and~\eqref{eqn8}, we get
\begin{equation*}
t_{a_4}^2t_{a_1}t_{a_2}t_{a_6}t_{a_8} t_{a_4}^2 t_{a_2}t_{a_3}t_{a_7}t_{a_8}
  =t_{a_3}t_{b_1}t_{b_2}t_{b_3}t_{b_7}t_{b_8}t_{b_6}t_{a_6}.
\end{equation*}
Cancelling $t_{a_3}$ and $t_{a_6}$ yields
\begin{equation*}
t_{a_4}^4 t_{a_1} t_{a_2}^2t_{a_7}t_{a_8}^2
    =t_{b_1}t_{b_2}t_{b_3}t_{b_7}t_{b_8}t_{b_6}.
\end{equation*}
Any simple closed curve on the left hand side is disjoint from
each closed curve on the right. Notice also that the complements
of $a_1\cup b_1$, $a_2\cup b_2$, $a_2\cup b_3$, $a_7\cup b_7$, $a_7\cup b_8$
and of $a_8\cup b_6$ are all connected. Lemma~\ref{Proposition} now implies
that $t_{a_4}^4$ is a product of three commutators, implying (b).
\end{proof}

\begin{prop}\label{negy}
Let $h\geq 2$ and let $a$ and $b$ be two simple closed curves intersecting
each other transversely at one point on $F$. Then $t_a^4t_b^4$ is a product
of three commutators.
\end{prop}
\begin{proof}
Suppose  that the two-holed torus of Figure~\ref{lantern}~(ii) is embedded
in $F$ in such a way that $a_4$ and $a_5$ are nonseparating on $F$. The curve
$a_2$ intersects $t_{a_1}(a_2)$ transversely at one point. Since $a$
intersects $b$ transversely at one point also and since any two such pairs
are topologically equivalent, we can assume that $a=a_2$ and $b=t_{a_1}(a_2)$.
By the two-holed torus relation, we have $t_{a_4}t_{a_5}=(t_{a_1}t_{a_2}t_{a_3})^4$.
 Let us denote $t_{a_i}$ by $t_i$. Then, we obtain
 \begin{eqnarray*}
 t_4t_5
 &=&    (t_1t_2t_3 t_1t_2t_3 )(t_1t_2t_3 t_1t_2t_3)\\
 &=&    (t_1t_2t_1 t_3t_2t_3 )(t_1t_2t_1 t_3t_2t_3) \\
 &=&    (t_2t_1t_2 t_2t_3t_2 )(t_2t_1t_2 t_2t_3t_2)\\
 &=&          t_2t_2t_3t_2t_2t_1t_2t_2t_3t_2t_2t_1 \\
 &=&    (t_2t_2t_3t_2^{-1}t_2^{-1}) t_2t_2 t_2t_2(t_1 t_2t_2t_2t_2t_1^{-1})
 t_1t_1(t_1^{-1}t_2^{-1}t_2^{-1}t_3t_2t_2t_1).
 \end{eqnarray*}
 If $v=t_{a_2}^2(a_3)$ and $w=t_{a_1}^{-1}t_{a_2}^{-2}(a_3)$, we have the
 equality
 \begin{eqnarray*}
 (t_{a_4} t_v^{-1} t_{a_5}t_{w}^{-1}) = t_a^4 t_b^4 t_{a_1}^2.
 \end{eqnarray*}
 Now, $t_{a_4} t_v^{-1} t_{a_5}t_{w}^{-1}$ is a commutator and
 $t_{a_1}^2$ is a product of two commutators. This observation
 completes the proof of Proposition~\ref{negy}.
\end{proof}

\section{Signature computations}
\label{fourth}
\input

The relations expressing certain products of Dehn twists
as products of commutators proved in Section~\ref{third}
allow us to construct corresponding Lefschetz fibrations.
These fibrations, and their signatures, depend on the
choices we make for the diffeomorpisms occurring in the
commutator relations.

In this section the genus of the fiber $F$ is $h\geq 3$.
The base $B$ of genus $g$ will be denoted by $\Sigma_g$
if it is closed, and by $\Sigma_g^r$ if it has $r$ boundary
components. For a smooth surface bundle $X\rightarrow\Sigma_g^r$,
the signature is completely determined by the corresponding
monodromy representation.
We shall pass back and forth between surface bundles over bases
with boundary and Lefschetz fibrations over closed bases using
the following well-known fact, see~\cite{GS,Mats,O}:
\begin{prop}\label{p:ss}
The signature of a fibered neighbourhood of a nonseparating,
respectively separating, singular fiber in a Lefschetz fibration
is equal to $0$, respectively to $-1$. \qed
\end{prop}

Now fix a symplectic basis for $H_1(F,\Z )$, so that the
monodromy representation
$$
\rho\colon\pi_1(\Sigma_g^r)\rightarrow\Gamma_h
$$
of $X$ composed with the action of the mapping class group on homology
$$
\phi\colon\Gamma_h\rightarrow Sp(2h,\Z )
$$
yields a symplectic representation $\chi$ of the fundamental group of the
base. The following result of Meyer~\cite{M} allows us to calculate the
signature:
\begin{thm}\label{t:Meyer}
Let $f\colon X\rightarrow\Sigma_g^r$ be an oriented surface bundle with
monodromy representation $\rho\colon\pi_1(\Sigma_g^r)\rightarrow\Gamma_h$.
Fix a standard presentation of $\pi_1(\Sigma_g^r)$ as follows:
$$
\pi_1(\Sigma_g^r)=\langle a_1,b_1,\ldots,a_g,b_g,c_1,\ldots,c_r \ \vert \
\prod_{i=1}^g[a_i,b_i]\prod_{j=1}^rc_j=1\rangle \ ,
$$
and let $\tau_h\colon Sp(2h,\Z )\times Sp(2h,\Z )\rightarrow \Z$ by
the cocycle defined in~\cite{M}.

Then the signature of $X$ is given by the formula
\begin{alignat*}{1}
\sigma (X)=\sum_{i=1}^g\tau_h(\kappa_i,\beta_i) &-
\sum_{i=2}^g\tau_h(\kappa_1\ldots\kappa_{i-1},\kappa_i) \\
&-\sum_{j=1}^{r-1}\tau_h(\kappa_1\ldots\kappa_g
\gamma_1\ldots\gamma_{j-1},\gamma_j) \ ,
\end{alignat*}
where $\alpha_i=\chi (a_i)$, $\beta_i=\chi (b_i)$,
$\gamma_i=\chi (c_i)$ and $\kappa_i=[\alpha_i,\beta_i]$.
\qed
\end{thm}

Here is a first application of this formula:

\begin{prop}\label{p:sDehn}
There is a Lefschetz fibration $X\rightarrow\Sigma _2$ with a unique
singular fiber and with signature $-1$, whether the vanishing
cycle is separating or not.
\end{prop}
\begin{proof}
It is well-known that a Dehn twist can be written as a product of
two commutators, see~\cite{KO}. We need to make explicit choices
for these commutators. To this end we consider curves $a$, $a_1$,
$a_2$, $a_3$, $b_1$, $b_2$ and $b_3$ on a
genus $3$ subsurface of $F$ as in Figure~\ref{KO}.
  \begin{figure}[hpbt]
  \centerline{\epsfbox{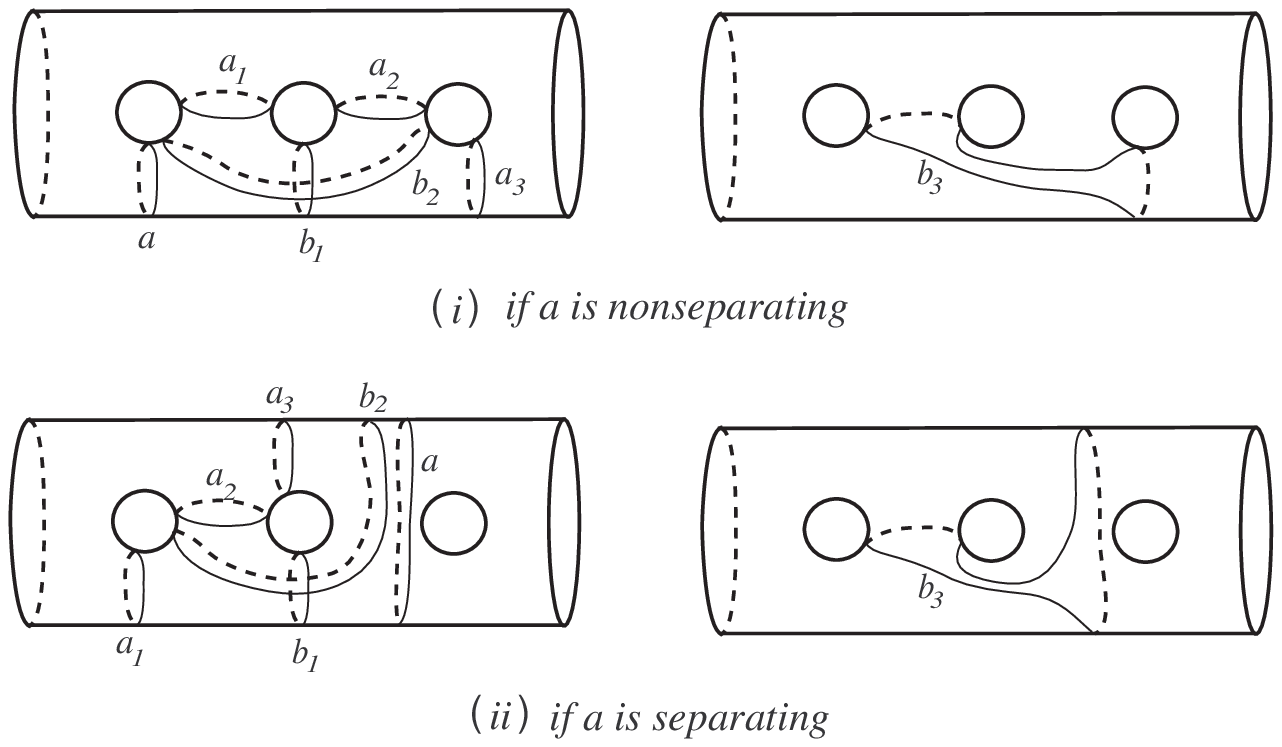}}
  \caption{}
  \label{KO}
  \end{figure}
Further, we add curves according to Figure~\ref{KOadd}.
  \begin{figure}[hpbt]
  \centerline{\epsfbox{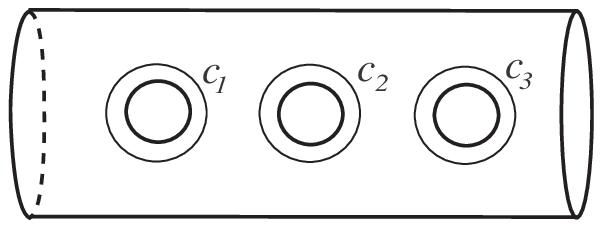}}
  \caption{}
  \label{KOadd}
  \end{figure}
If we choose the genus $3$ subsurface suitably, the vanishing cycle $v$
is topologically equivalent to the curve $a$.

Define diffeomorphisms $\phi_1$ and $\phi_2$ of $F$ as follows. If
$a$ is nonseparating, set
 $$
 \phi_1=t_{c_1}t_{b_2}t_{c_2}t_{a_2}t_{b_1}t_{c_2}t_{a_1}t_{c_1}
 $$
and
 $$
  \phi_2 = t_{c_3}t_{b_3}t_{a_3}t_{c_3} \ .
 $$ If $a$ is separating, set
 $$
 \phi_1=t_{c_2}t_{a_2}t_{c_1}t_{b_2}t_{a_1}t_{c_1}t_{b_1}t_{c_2}
$$ and $$ \phi_2= t_{c_2}t_{a_3}t_{b_3}t_{c_2} \ . $$

 One can check that $\phi_1(a_1)=b_2$, $\phi_1(b_1)=a_2$ and
$\phi_2(a_3)=b_3$.

The lantern relation as in Theorem 2 of~\cite{KO} implies that
 $$t_{a_3}t_{b_3}^{-1}t_{a_2}t_{b_2}^{-1}t_{a_1}t_{b_1}^{-1}t_{a}=1
 \ .$$

The monodromy representation of the complement of the singular
fiber is given by mapping the standard generators of
$\pi_1(\Sigma_2^1)$ to $t_{a_3}$, $\phi_2$, $\phi_1$,
$t_{a_1}^{-1}t_{b_1}$ and $t_a$ respectively, as
 $$
 [t_{a_3},\phi_2][\phi_1,t_{a_1}^{-1}t_{b_1}]t_{a}=1 \ .
 $$
Evaluating the signature cocycle, Theorem~\ref{t:Meyer} shows that
the complement of the singular fiber has signature $-1$ if $v$ is
nonseparating, and has signature $0$ if $v$ is separating. Now
Proposition~\ref{p:ss} and Novikov additivity complete the proof.
\end{proof}

{\it Mutatis mutandis}, this calculation generalizes to prove the
next three Propositions:
\begin{prop}\label{p:sDehn^2}
There is a Lefschetz fibration $X\rightarrow\Sigma _2$ with two
singular fibers whose monodromies are Dehn twists with the same
nonseparating vanishing cycle and signature equal to $-2$.
\end{prop}

\begin{prop}\label{p:sDehn^4}
There is a Lefschetz fibration $X\rightarrow\Sigma _3$ with four
singular fibers whose monodromies are Dehn twists with the same
nonseparating vanishing cycle and signature equal to $-4$.
\end{prop}

\begin{prop}\label{p:sta^4tb^4}
Let $a$ and $b$ be two nonseparating simple closed curves on $F$
which intersect transversely and precisely at one point. There is
a Lefschetz fibration $X\rightarrow\Sigma_3$ with signature $-4$
which has eight singular fibers, four of which have monodromy a
Dehn twist along $a$ and four of which have monodromy a Dehn twist
along $b$.
\end{prop}

\begin{proof}[Proof of Propositions~\ref{p:sDehn^2}-\ref{p:sta^4tb^4}]
In all these proofs the signature of the Lefschetz fibration
is the same as that of the complement of the singular fibers,
because all the vanishing cycles are nonseparating,
cf.~Proposition~\ref{p:ss}.

We take curves $x$, $y$, $a_1,\ldots,a_8$, $b_1,\ldots,b_8$ on a
genus $3$ subsurface of $F$ as in Figure~\ref{feners}.
We also add curves $c_1$, $c_2$, $c_3$ as in Figure~\ref{KOadd},
and $d$ and $e$ as in Figure~\ref{commadd}.
\begin{figure}[hpbt]
\centerline{\epsfbox{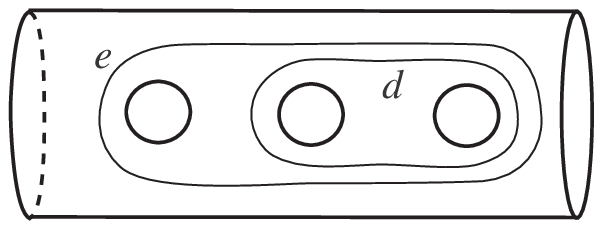}}
\caption{}
\label{commadd}
\end{figure}

For each of the Propositions, the vanishing cycles $a$ and $b$ are
topologically equivalent to certain curves $a_0$ and $b_0$. We fix
the latter explicitly, and construct some  diffeomorphisms as
required by the proofs in Section~\ref{third}, so that we can
write the monodromy representation of the complement of the
singular fibers as a relator in the mapping class group. Then the
calculation is done by implementing the formula in Theorem~\ref{t:Meyer}
with the following data.

For Proposition~\ref{p:sDehn^2} the base is $\Sigma_2^2$ and $a_0$
is taken to be $a_4$. The relator giving the monodromy
representation is
$$
[t_{b_3}^{-1}t_{a_1},\phi_1][t_{b_1}^{-1}t_{a_2},\phi_2]t_{a_4}^2=1
 \ ,
$$
with
$$
\phi_1=t_{e}t_{b_2}t_{a_1}t_{e}t_{c_2}t_{a_6}t_{b_3}t_{c_2}
$$
and
$$
\phi_2=t_{c_1}t_{a_3}t_{a_2}t_{c_1}t_{d}t_{a_8}t_{b_1}t_{d} \ .
$$

For Proposition~\ref{p:sDehn^4} the base is $\Sigma_3^4$ and $a_0$
is taken to be $a_4$ again. The relator giving the monodromy
representation is
$$
 [t_{b_6}^{-1}t_{a_2},\phi_1][t_{b_7}^{-1}t_{a_8},\phi_2]
 [t_{b_2}^{-1}t_{a_1},\phi_3]t_{a_4}^4=1
$$
with
$$
\phi_1=t_{c_1}t_{b_8}t_{a_2}t_{c_1}t_{c_3}t_{a_8}t_{b_6}t_{c_3} \ ,
$$
$$
\phi_2=t_{c_3}t_{b_3}t_{a_8}t_{c_3}t_{c_1}t_{a_2}t_{b_7}t_{c_1}
$$
and
$$
\phi_3=t_{c_1}t_{b_1}t_{a_1}t_{c_1}t_{c_3}t_{a_7}t_{b_2}t_{c_3} \ .
$$

For Proposition~\ref{p:sta^4tb^4} consider the curves
$a_1,\ldots,a_6$ on a genus $2$ subsurface of the fiber as in
Figure~\ref{lantern}~(ii). 
The base surface is $\Sigma_3^8$ and the curve $a_0$ is taken to
be $a_2$ and $b_0$ is taken to be $t_{a_1}(a_2)$. We first compute
the signature corresponding to the relator $$
[t_{a_1}^{-1}t_{a_2}^{-2}t_{a_3}t_{a_2}^2t_{a_1}t_{a_5}^{-1},\phi
] t_{a_2}^4(t_{a_1}t_{a_2}t_{a_1}^{-1})^4t_{a_1}^2=1 \ , $$ with
$$ \phi =
t_{a_2}^3t_{a_3}t_{a_6}t_{a_4}t_{a_5}t_{a_6}t_{a_3}t_{a_2}
t_{a_6}t_{a_5}t_{a_3}t_{a_6}t_{a_2}^2t_{a_1} \ , $$ which is the
monodromy of a fibration $X'\rightarrow T^2$ with $10$ singular
fibers. The signature of $X'$ is equal to $-6$. Subtracting off
the fibration of Proposition~\ref{p:sDehn^2} from $X'$ gives the
claim.
\end{proof}

As the expression of a given element in $\Gamma_h$ as a product of
commutators is not unique, it is conceivable that the signature
of the corresponding Lefschetz fibrations might be different for
different choices of commutators. Then a surface bundle of nonzero
signature could be constructed by subtracting the Lefschetz fibrations
corresponding to different choices from each other.

\section{Bounds on the genus function $g_h(n)$}

We now prove the theorems about the minimal genus function
$g_h(n)$ stated in the Introduction.

\begin{proof}[Proof of Theorem~\ref{main1}]
We apply the subtraction operation to the Lefschetz fibrations
$X_1\to \Sigma _3$ and $X_2\to \Sigma _3$
as in Proposition~\ref{p:sta^4tb^4} and~\ref{p:sDehn^4},
respectively. In $X_1$ we group the singular fibers into two groups
each containing four singular fibers with coinciding vanishing cycles;
in $X_2$ the singular fibers form one group. Now subtracting two
copies of $X_2$ according to the above pattern we get a surface
bundle $Y_h\to \Sigma _9$ of fiber genus $h$ with
$\sigma (Y_h)=\sigma (X_1)-2\sigma (X_2)= -4-2(-4)=4$
(cf.~Propositions~\ref{p:sDehn^4} and \ref{p:sta^4tb^4}).
Thus $g_h(1)\leq 9$, and the claim now follows by pulling $Y_h$ back
to unramified coverings of $\Sigma _9$ of degree $\vert n\vert$.
\end{proof}
Surface bundles over $\Sigma _9$ with higher signature
can be constructed as follows.
\begin{proof}[Proof of Theorem~\ref{main2}]
Notice that the relators defining the fibrations we used in the
proof of Theorem~\ref{main1} represent 1 in the mapping class group
$\Gamma^1_h$ of a surface with one boundary component.
According to Proposition~\ref{secc}, this fact shows that  the
fibrations given by the relators
$\Pi _{i=1}^3 [ a_i, b_i ]t_a^4 t_b^4$ and $\Pi _{i=1}^3 [ c_i, d_i ]t_a^4$
admit sections with vanishing self-intersection.
Since the lifts of the various Dehn twists are chosen to be Dehn twists
in $\Gamma^1_h$, Proposition~\ref{osszerag} implies that
$Y_h\to \Sigma _9$ also
admits a section with zero self-intersection for all $h$.
Now write $h$ as $3k+l$ where $l\in \{ 0, 1,2 \}$, and apply
Lemma~\ref{secrag} to $k$ copies of $Y_3$
together with the product $\Sigma _l \times \Sigma _9
\to \Sigma _9$. The resulting surface bundle
$S_h\to \Sigma _9$ of fiber genus $h$
has  $\sigma (S_h)=4k=4\frac{h-l}{3}\geq 4\frac{h-2}{3}$.
\end{proof}

Now we  turn to the study of the
asymptotic behaviour of the genus function.

\begin{proof}[Proof of Theorem~\ref{main3}]
First notice that the proof of Theorem~\ref{main2} immediately yields
$G_h =\lim _{n\to \infty }\frac{g_h(n)}{n}\leq \frac{24}{h-l}$ for all
$h$. (As before,  $l\in \{ 0,1,2\} $ is the mod 3 residue of $h$.)

Now every surface of odd genus is a covering of a genus $3$
surface. It was shown in Lemma 4.1 of~\cite{Morita}, that after
replacing a given surface bundle by a pullback to some covering of the
base, the resulting surface bundle admits fiberwise coverings of
any given degree. From the multiplicativity of the signature
in coverings and the multiplicativity of the Euler characteristic
of the fiber in fiberwise coverings, for odd $h$ we obtain
\begin{equation}\label{fiberwise}
\lim_{n\to\infty}\frac{g_h(n)}{n}\leq\frac{2}{h-1}
\lim_{n\to\infty}\frac{g_3(n)}{n} \leq \frac{16}{h-1} \ .
\end{equation}

For even $h$ consider the fibration $Z \to \Sigma _{8n+1}$
of fiber genus $h-1$  with signature $2n(h-2)$ we got by taking
fiberwise coverings.
It is easy to see that since $Y_3\to \Sigma _9$ admits a section
of zero self-intersection, so does $Z\to \Sigma _{8n+1}$. Summing $Z$ and
the product fibration $\Sigma _1\times \Sigma _{8n+1}\to \Sigma _{8n+1}$
along their sections (as in Lemma~\ref{secrag}), we get a
fibration over $\Sigma _{8n+1}$ with fiber genus $h$
and signature $2n(h-2)$. These examples yield the bound
$G_h\leq \frac{16}{h-2}$ once $h$ is even.
Consequently the proof of Theorem~\ref{main3} is complete.
\end{proof}

\begin{rem}
For certain values of $h$, the examples of Kodaira~\cite{Kd}
give a better upper bound, namely $G_h\leq\frac{44}{5(h-1)}$.
Our construction has the advantage of covering all possible
values of $h$ (and $n$).
\end{rem}


\begin{thebibliography}{AAAA}

\bibitem{A}
M.~F.~Atiyah,
{\em The signature of fibre-bundles},
Global Analysis, Papers in Honor of K.~Kodaira, Tokyo Univ.~Press,
1969, 73--84.

\bibitem{CHS}
S.~S.~Chern, F.~Hirzebruch and J.~P.~Serre,
{\em On index of a fibered manifold},
Proc.~Amer.~Math.~Soc.~{\bf8} (1957), 587--596.

\bibitem{E}
H.~Endo,
{\em A construction of surface bundles over surfaces with non-zero
signature},
Osaka J.~Math.~{\bf 35} (1998), 915--930.

\bibitem{GS}
R.~E.~Gompf and A.~I.~Stipsicz,
{\em  4-manifolds and Kirby calculus},
Graduate Studies in Mathematics, vol. {\bf20}, American Math.~Society,
Providence 1999.

\bibitem{H}
J.~Harer, {\em The second homology group of the mapping class group
of an orientable surface}, Ivent.~Math.~{\bf 72} (1983), 221--239.

\bibitem{Ho}
M.~Hoster, {\em A new proof of the signature formula for surface
bundles}, Top.~and~App.,~to appear.

\bibitem{J} D.~L.~Johnson, {\em Homeomorphisms of a surface which act
trivially on homology,} Proc.~Amer.~Math.~Soc.~{\bf75}  (1979),
119-125.

\bibitem{Ki}
R.~Kirby,
{\em Problems in low-dimensional topology}, in Geometric Topology
(W. Kazez ed.) AMS/IP Stud.~Adv.~Math.~vol 2.2, American Math.~Society,
Providence 1997.

\bibitem{Kd}
K.~Kodaira,
{\em A certain type of irregular algebraic surfaces},
J.~Anal.~Math. {\bf19} (1967), 207--215.

\bibitem{KO}
M.~Korkmaz and B.~Ozbagci, {\em Minimal number of singular fibers in a
Lefschetz fibration}, Proc.~AMS {\bf129} (2001), 1545--1549.

\bibitem{Kt}
D.~Kotschick,
{\em Signatures, monopoles and mapping class groups},
Math.~Res.~Letters {\bf5} (1998), 227--234.

\bibitem{La}
K.~Lamotke,
{\em The topology of complex projective varieties after S.~Lefschetz},
Topology {\bf 20} (1981), 15--51.

\bibitem{Li}
W.~B.~R.~Lickorish,
{\em A finite set of generators for the homeotopy group of a $2$-manifold},
Math.~Proc.~Camb.~Phil.~Soc.~{\bf60} (1964), 769--778.

\bibitem{Mats}
Y.~Matsumoto, {\em Lefschetz fibrations of genus two --
a topological approach},
Proc.~of the 37th Taniguchi symposium on ``Topology and Teichm\"uller spaces'',
World Scientific 1996.

\bibitem{M}
W.~Meyer,
{\em Die Signatur von Fl\"achenb\"undeln}, Math.~Ann.~{\bf201} (1973),
239--264.

\bibitem{Morita}
S.~Morita,
{\em Characteristic classes of surface bundles}, Invent.~Math.~{\bf 90}
(1987), 551--577.

\bibitem{O}
B.~Ozbagci,
{\em Signatures of Lefschetz fibrations},
Pacific J. of Math., to appear.

\bibitem{Sm}
I.~Smith,
{\em Symplectic geometry of Lefschetz fibrations},
D.~Phil.~thesis, Oxford 1998.

\bibitem{S}
A.~I.~Stipsicz,
{\em Surface bundles with nonvanishing signatures},
Preprint.

\bibitem{Waj}
B.~Wajnryb,
{\em An elementary approach to the mapping class group of a surface},
Geometry and Topology {\bf3} (1999), 405--466.

\end{thebibliography}
\end{document}